\documentclass[12pt]{amsart}

\usepackage{amsthm,amsmath,amssymb}
\usepackage{appendix}
\usepackage{enumerate}

\theoremstyle{plain}
\newtheorem{theorem}{Theorem}
\newtheorem{lemma}{Lemma}
\newtheorem{corollary}{Corollary}
\newtheorem{proposition}{Proposition}
\newtheorem{conjecture}{Conjecture}

\theoremstyle{definition}

\newtheorem{example}{Example}

\newcommand{\raf}[1]{(\ref{#1})}

\newcommand{\cG}{{\mathcal G}}
\newcommand{\cH}{{\mathcal H}}

\newcommand{\cF}{{\mathcal{F}}}
\newcommand{\cP}{{\mathcal P}}

\newcommand{\ga}{\alpha}

\def\G{\mathcal{G}}

\def\T{\mathcal{T}}
\def\ZZ{\mathbb{Z}}
\def\Zp{\mathbb{Z}_{+}}
\def\mex{\operatorname{mex}}
\def\NIM{{{\sc Nim}}}

\begin{document}

\title{On the Sprague-Grundy Function of Tetris Extensions of Proper {\sc Nim}}

\author{Endre Boros}
\address{MSIS and RUTCOR, RBS, Rutgers University,
100 Rockafeller Road, Piscataway, NJ 08854}
\email{endre.boros@rutgers.edu}

\author{Vladimir Gurvich}
\address{MSIS and RUTCOR, RBS, Rutgers University,
100 Rockafeller Road, Piscataway, NJ 08854;
National Research University Higher School of Economics (HSE), Moscow}
\email{vladimir.gurvich@rutgers.edu}

\author{Nhan Bao Ho}
\address{Department of Mathematics and Statistics, La Trobe University, Melbourne, Australia 3086}
\email{nhan.ho@latrobe.edu.au, nhanbaoho@gmail.com}

\author{Kazuhisa Makino}
\address{Research Institute for Mathematical Sciences (RIMS)
Kyoto University, Kyoto 606-8502, Japan}
\email{e-mail:makino@kurims.kyoto-u.ac.jp}



\begin{abstract}
Given a hypergraph  $\cH \subseteq 2^I \setminus \{\emptyset\}$ on the ground set
$I = \{1, \ldots, n\}$, we assign to each  $i \in I$
a nonnegative integer  $x_i$, that is a pile of $x_i$ tokens, and
consider the following generalization of the classical game of  {\sc Nim}:
Two players alternate turns. In a move a player chooses an arbitrary edge
$H \in \cH$  and reduces all piles $i \in H$.
The player who is out of moves loses.
We call the obtained game hypergraph  {\sc Nim}.
Such a game is called proper {\sc Nim}, when
$\cH=2^I \setminus\{I,\emptyset\}$ is the family of all proper subsets of $I$.
Jenkyns and Mayberry \cite{JM80} described the Sprague-Grundy
(or SG in short) function of these games.
In this paper we introduce Tetris extensions of hypergraph {\sc Nim}, and
obtain a closed formula for the SG functions of the extensions of proper {\sc Nim}, when
$n\geq 3$. Surprisingly, the case of $n=2$ is much more complicated.
For this case we only suggest several partial results and conjectures.
\end{abstract}

\maketitle

\section{Introduction}
\label{s0}

In the classical game of \NIM~ there are
$n$ piles of stones and two players move alternating. A move consists of
choosing a nonempty pile and taking some positive number of stones from it.
The player who cannot move is the looser.
Bouton \cite{Bou901} analyzed this game and described the winning strategy for it.

In this paper we consider the following generalization of \NIM.
Given a hypergraph $\cH \subseteq 2^I$, where $I=\{1 \ldots ,n\}$
and a \emph{position} $x = \{x_1, \ldots, x_n\} \in \ZZ_+^I$,
two players alternate in choosing a hyperedge
$H\in \cH$ and strictly decreasing values of $x_i$  for all  $i\in H$.
We assume in this paper that $\emptyset\not\in \cH$ for all  considered hypergraphs.
The player who cannot move is loosing.
We call the obtained game \emph{hypergraph} {\sc Nim} and denote it by \NIM$_\cH$.

These games belong to the class of impartial games.
In this paper we do not need to immerse in the theory of impartial games.
We will need only to recall a few basic facts to explain and motivate our research.
We refer the reader to \cite{Alb07,BCG01-04} for more details.

It is known that the set of positions of an impartial game can
uniquely be partitioned into sets of \emph{winning} and \emph{loosing} positions.
Every move from a loosing position goes to a winning one, while
from a winning position always there is a move to a loosing one.
This partition shows how to win the game whenever possible.
The so-called Sprague-Grundy (SG) function $\G_\Gamma$ of an impartial game $\Gamma$ provides
a refinement of the above partition.
Namely, for $x\in\ZZ_+^I$ we have $\G_\Gamma(x)=0$ if and only if $x$ is a loosing position.
The notion of the SG function for impartial games was introduced by Sprague and Grundy
\cite{Spr35,Spr37,Gru39} and it plays a fundamental role in the analysis of composite impartial games.

Finding a formula for the SG function of an impartial game remains a challenge, in general.
Closed form descriptions are known only for some special classes.
The purpose of our research is to extend these results and to describe classes of hypergraphs
for which we can provide a closed form for the SG function of \NIM$_\cH$.
To follow our proofs, we need to recall the precise definition of the SG function,
which we will do in the beginning of Section \ref{s2}.

Hypergraph \NIM~ generalizes several families of impartial games considered in the literature.
For instance, if $\cH=\{\{1\}, \ldots, \{n\}\}$, then
\NIM$_\cH$ is the classical \NIM, which was analyzed and solved by Bouton \cite{Bou901}.
The case of $\cH=\{S\subseteq I\mid 1 \leq |S| \leq k\}$, where $k<n$, was considered by
Moore \cite{Moo910}.
He characterized for these games the set of loosing positions.
Jenkyns and Mayberry \cite{JM80} described also the
set of positions in which the SG value is $1$, and
provided an explicit formula for the SG function in the case of $k = n-1$,
that is, for 
\emph{proper} \NIM.
Let us note that it seems difficult to extend this result for any $k$  such that $1<k<n-1$.
For instance, no closed formula is known for the SG function when $n=4$ and $k=2$.

In \cite{BGHMM15} the game \NIM$_\cH$ was considered in the case of $\cH=\{S\subseteq V\mid |S|=k\}$,
and the corresponding SG function was determined by a closed formula whenever  $2k\geq n$.
Let us add that even winning positions are not known when $2k < n$, e.g., for $n=5$ and $k=2$.

In this paper first we consider a simple extension operation of a hypergraph.
Let us set $I_+ = \{0,1,...,n\}=\{0\}\cup I$.
To a hypergraph $\cH\subseteq 2^I\setminus\{\emptyset\}$ we associate an extended hypergraph $\cH_+\subseteq 2^{I_+}$ defined by
\[
\cH_+ ~=~ \cH\cup \{0\}\cup \{H\cup\{0\}\mid H\in\cH\}.
\]
We call \NIM$_{\cH_+}$ the \emph{extension} of \NIM$_\cH$.
For clarity, we shall denote by $(x_0;x)$ the positions of the extended game, where $x\in\ZZ_+^I$.

Even though this is a very simple operation,
the SG function of such an extension can be much more complicated than the SG function of the original game.
See for instance, the extension of $2$-pile \NIM, studied in Section \ref{s3}.

To state our first main result we need a few more definitions. To a position $x\in\ZZ_+^I$ we associate
\begin{equation}\label{e-mu}
m(x)=\min_{i\in I} x_i ~~~\text{and}~~~ u(x)=\sum_{i\in I} x_i.
\end{equation}
Furthermore, we say that a hypergraph $\cH\subseteq 2^I$ satisfies the \emph{min-sum} property, if for the SG function $\G_\cH$ of \NIM$_\cH$ we have
\begin{equation}\label{e-minsum}
\G_\cH(x) ~=~ F(m(x),u(x))
\end{equation}
for some function $F:\ZZ_+^2\to \ZZ_+$.
We call $\cH$ \emph{hereditary} if for all $H\in\cH$ and $\emptyset\neq H'\subseteq H$ we have $H'\in\cH$.

\begin{theorem}\label{t1}
Assume that $n\geq 3$ and $\cH\subseteq 2^I$ is a hereditary hypergraph satisfying the min-sum property \eqref{e-minsum}.
Assume further that all inclusion-wise maximal hyperedges of $\cH$ are of size at least $2$.
Then for the SG function of the extended \NIM$_{\cH_+}$ game we have
\[
\G_{\cH_+}(x_0;x) ~=~ F(m(x),u(x)+x_0).
\]
\end{theorem}

We can apply this result to a proper \NIM~ game with $n\geq 3$ and
obtain a closed formula for the SG function of its extension.

To a position $(x_0;x)\in\ZZ_+^{I_+}$ let us associate
\begin{equation}\label{e-yz}
y(x_0;x)=u(x)+x_0-n\cdot m(x) ~~~\text{and}~~~ z(x_0;x)=\binom{y(x_0;x)+1}{2}+1.
\end{equation}

\begin{corollary}\label{c1}
Let $n\geq 3$ and $\cH=2^I\setminus\{I,\emptyset\}$.
Then for the extension of the proper \NIM~ game defined by this hypergraph, we have
\begin{equation}\label{e-uv}
\G_{\cH_+}(x_0;x) = \begin{cases}
u+x_0 &\text{if~} m<z,\\
(z-1)+((m-z) \mod (y+1) &\text{if~} m\geq z,
\end{cases}
\end{equation}
where for simplicity we use $m=m(x)$, $u=u(x)$, $y=y(x_0;x)$ and $z=z(x_0;x)$.
\end{corollary}

We will illustrate this statement by Example \ref{example1} in Section \ref{s2}.

Note that when $x_0=0$, the above formula coincides
with the SG function of proper \NIM, obtained by Jenkyns and Mayberry \cite{JM80}.
Note also that $\cH$ satisfies all conditions of Theorem \ref{t1}.
Thus, the above claim follows from the result of \cite{JM80} and Theorem \ref{t1}.

Let us add that for any $n\geq 2$ the set

\medskip

$\cP(\cH_+)=\{(x_0;x)\mid x_0=0, \mbox{ and } x_1=\cdots =x_n\}$

\smallskip
\noindent
is the set of loosing positions of \NIM$_{\cH_+}$.
This can be easily seen from the characterization of loosing positions given above.
For $n \geq 3$ one can also verify that
$\G_{\cH_+}(x_0;x)=0$  if and only if $(x_0;x)\in \cP(\cH_+)$.

\smallskip

We generalize Theorem \ref{t1} further, by observing that
the single $x_0$ component in the extension can be replaced by a hypergraph \NIM~ game
that satisfies certain properties.
For this let us consider a finite set $J$, disjoint from $I$ and a hypergraph $\cF\subseteq 2^J$.
We define the extension of $\cH\subseteq 2^I$ by $\cF$ as follows:
\[
\cH \otimes \cF ~=~ \cH \cup \cF \cup \cH\times\cF,
\]
where $\cH\times\cF=\{H\cup F\mid H\in\cH,~ F\in\cF\}$ is the usual direct product.
If $J=\{0\}$ and $\cF=\{\{0\}\}$, then the above extension coincides with $\cH_+$.

For \NIM$_\cF$ and a position $x^J\in \ZZ_+^J$ we denote by $\T_\cF(x^J)$ the maximum number of consecutive moves the players could make in \NIM$_\cF$ starting form $x^J$.
This function was introduced and studied in \cite{BGHMM15,BGHMM16}, and called the \emph{Tetris function} defined by the hypergraph $\cF$. Furthermore, we call the hypergraph $\cF$ a \emph{Tetris hypergraph} if the SG function of \NIM$_\cF$ coincides with $\T_\cF$.

\begin{theorem}\label{tt2}
Assume that $n\geq 3$, $\cH\subseteq 2^I$ is a hereditary hypergraph satisfying
the min-sum property \eqref{e-minsum}, and
all inclusion-wise maximal hyperedges of $\cH$ are of size at least $2$.
Assume further that $\cF\subseteq 2^J$ ($I\cap J=\emptyset$) is a Tetris hypergraph.
Then for the SG function of the extended \NIM$_{\cH\otimes\cF}$ game we have
\[
\G_{\cH\otimes\cF}(x^J;x^I) ~=~ F(m(x^I),u(x^I)+\T_\cF(x^J)).
\]
\end{theorem}
We call \NIM$_{\cH\otimes\cF}$ the \emph{Tetris extension} of \NIM$_\cH$ if $\cF$ is a Tetris hypergraph.
For instance, $\cF=\{\{0\}\}$ corresponding to a single pile \NIM~ is Tetris. Thus, \NIM$_{\cH_+}$ is a Tetris extension.
Another example for Tetris hypergraphs is an intersecting hypergraph, that is such that
$F\cap F'\neq\emptyset$ for any two hyperedges $F,F'\in\cF$.
More examples are described in \cite{BGHMM16}.

We can extend a proper \NIM~ game with any of these Tetris hypergraphs and
obtain a closed formula for the SG function of those extensions
by the result of Jenkyns and Mayberry \cite{JM80} and Theorem \ref{tt2}.

\bigskip

The rest of the paper is organized as follows.
In Section \ref{s2} we prove Theorems \ref{t1} and \ref{tt2}.
In Section \ref{s3} we study the extension of proper \NIM~ with two piles.
In this case the SG function behaves in a chaotic way, and
we only obtain some partial results and state some conjectures.
Finally, in Section \ref{s4} we consider the SG function of
a game {\em slow \NIM} in which the size of a pile can be reduced by at most one.


\section{Tetris extensions of proper \NIM~ for $n\geq 3$}
\label{s2}
In this section we prove Theorems \ref{t1} and \ref{tt2}.
For this let us recall first the definition of the SG function of an impartial game.
For a subset $S\subseteq \ZZ_+$ of nonnegative integers
we associate the value of the smallest integers, not belonging to $S$:
\[
mex(S) ~=~ \min_{v\in \ZZ_+\setminus S} v.
\]
In particular, $mex(\emptyset)=0$.

To an impartial game $\Gamma$ we associate its Sprague-Grundy function
$\G_\Gamma$ that assigns to every position $x$
of the game a nonnegative integer defined by the following recursion
\[
\G_\Gamma(x) ~=~ mex \{ \G_\Gamma (x') \mid \forall x' \text{ such that there is a move } x \to x'\}.
\]
Equivalently, the SG function can be defined by the following two properties:
\begin{itemize}
\item[(i)] For any move $x\to x'$ we have $g(x')\neq g(x)$.
\item[(ii)] For any position $x$ with $g(x)>0$ and for any integer $0\leq v< g(x)$ there is a move $x\to x'$ such that $g(x')=v$.
\end{itemize}
Any nonnegative integer valued function $g$ satisfying the above properties must coincide with $\G_\Gamma$.

Let us consider \NIM$_\cH$ and its extension, as in Theorem \ref{t1}.
For a position $(x_0;x)\in\ZZ_+^{I_+}$, let $k$ be an index such that $x_k=m(x)$ and $1 \leq k \leq n$, where $m(x)$ is defined in \eqref{e-mu}.
To such a position we associate a set of positions of \NIM$_\cH$ as follows
\begin{equation}\label{e-Z}
Z(x_0;x)=\{ z\in\ZZ_+^I \mid z_k=x_k,~ z_i \geq x_i ~(\forall i\not=k), ~~  u(z)=u(x)+x_0\}.
\end{equation}
By the above definition, it holds that
\begin{equation}
\label{eq-new1a}
m(z)=m(x)=x_k \mbox{ and } u(z)=u(x)+x_0 \mbox{ for all } z \in Z(x_0;x).
\end{equation}
Let us further define
\begin{equation}\label{e-A}
A(x_0;x)=\{(m(z'),u(z')) \mid \exists z \in Z(x_0;x), \exists \mbox{ a move } z \to z' \mbox{ in \NIM}_{\cH}\}.
\end{equation}

\begin{lemma}\label{lemma-1az-1}
Assume that $\cH$ satisfies the min-sum property \eqref{e-minsum}.
Then, for any position $(x_0;x)$ of the extension \NIM$_{\cH_+}$,
we have
\begin{itemize}
\item[{\rm (I)}] $F(m(x), u(x)+x_0) \not\in F(A(x_0;x))$, and
\item[{\rm (II)}] $\{0,1, \dots ,  F(m(x), u(x)+x_0)-1\} \subseteq F(A(x_0;x))$,
\end{itemize}
where $F(A(x_0;x))=\{F(m',u')\mid (m',u')\in A(x_0;x)\}$.
\end{lemma}

\proof
This follows from \raf{eq-new1a} and the min-sum property of $\cH$.
\qed

\begin{lemma}\label{lemma-1az0}
Assume that $\cH$ is a hereditary hypergraph, and $(x_0;x)\in\ZZ_+^{I_+}$ is a position of the extension \NIM$_{\cH_+}$.
For any $z \in Z(x_0;x)$ and any move $z\to z'$ in \NIM$_{\cH}$,
there exists a move $(x_0;x)\to (x_0';x')$ in the extension such that $m(x')=m(z')$ and  $u(x')+x'_0=u(z')$.
\end{lemma}

\proof
Assume that $k\in I$ is the index we used in the definition of $Z(x_0;x)$.
Let $z_i=x_i+\ga_i$ for $i\in I$ such that $\ga_k=0$, $\sum_{i\in I}\ga_i=x_0$, and $\ga_i\geq 0$ for all $i \in I$.
Let us define $x'\in\ZZ_+^I$ by
\[
x_i' = \begin{cases}
x_i & \mbox{if } z_i' \geq x_i,\\
z_i' & \mbox{otherwise},
\end{cases}
\]
and set $x_0'=\sum_{i\in I}\max \{z_i'-x_i, 0\}$.
Then $(x_0';x')\in\ZZ_+^{I_+}$, $x' \leq x$ and $x_0'\leq x_0$.
Define $S=\{i\in I\mid x'_i<x_i\}$, and note that $S$ is contained in the hyperedge $\{i\in I\mid z'_i<z_i\}$ of $\cH$ by the above construction.
Thus, if $S\neq\emptyset$ then $(x_0;x)\to (x_0';x')$ is a move in \NIM$_{\cH_+}$ by the hereditary property of $\cH$.
If $S=\emptyset$, then $x'_0<x_0$ because the move $z\to z'$ reduces at least one component strictly,
and thus $(x_0;x)\to (x_0';x')$ is also a move in \NIM$_{\cH_+}$.
\qed

\begin{lemma}\label{lemma-1az1}
Assume that $n\geq 3$, $\cH\subseteq 2^I$ is hereditary, and any element $i \in I$ is contained in a hyperedge of size at least $2$.
Let us consider a position $(x_0;x)\in\ZZ_+^{I_+}$ of \NIM$_{\cH_+}$.
For any  move $(x_0;x)\to (x_0';x')$ in \NIM$_{\cH_+}$, there exist a $z \in Z(x_0;x)$ and a move $z\to z'$ in \NIM$_{\cH}$ such that
$z' \in Z(x_0';x')$.
\end{lemma}

\proof
Assume that $k\in I$ is the index we used in the definition of $Z(x_0;x)$.
Furthermore, let $h$ denote an index such that $x'_h=m(x')$ and $h\in I$.
We consider separately four cases.

\smallskip

Case 1: $x_0'< x_0$ and $x_i' = x_i$ for all $i \in I$.
For an arbitrary index  $j\in I$,
we define positions $z,z'\in\ZZ_+^I$ by
\begin{equation}\label{e-zzp}
z_i = \begin{cases}
x_j+x_0 & \mbox{if } i=j,\\
x_i & \mbox{otherwise},
\end{cases}
~~~~
z_i' = \begin{cases}
x'_j+x'_0 & \mbox{if } i=j,\\
x'_i & \mbox{otherwise}.
\end{cases}
\end{equation}
Note that $\{j\} \in \cH$ by the assumption on $\cH$.
Since $z'_j < z_j$ and $z_i' = z_i$ for all $i \in I \setminus \{j\}$,
we have $z \in Z(x_0;x)$ and $z\to z'$ is a move in \NIM$_{\cH}$.

\smallskip

Case 2: There is an index $j \in I$ such that $j\not= k,h$ and $x_j' < x_j$.
In this case, we again can use positions $z,z'\in\ZZ_+^I$, defined in \eqref{e-zzp},
since
\begin{equation}\label{e-zpzxpx}
\{i \in I \mid z'_i<z_i\} = \{i \in I \mid x'_i<x_i\} \in \cH,
\end{equation}
implies that $z \in Z(x_0;x)$ and $z\to z'$ is a move in \NIM$_{\cH}$.

\smallskip

Case 3: $k=h$, $x'_k<x_k$ and $x'_i=x_i$ for all $i \in I\setminus \{k\}$.
By the assumption that any element is contained in hyperedge of size at least $2$ and since $\cH$ is hereditary,
there must exist an index $j \in I\setminus \{k\}$ such that $\{k,j\} \in \cH$.
Similarly to the previous case we consider the positions $z,z'\in\ZZ_+^I$ as defined in \eqref{e-zzp}.
Then we have again $z \in Z(x_0;x)$ and that $z\to z'$ is a move in \NIM$_{\cH}$, since
$\{i\mid z'_i<z_i\}\subseteq \{k,j\}$ and $\cH$ is hereditary.

\smallskip

Case 4: $k \not=h$ and $x'_i=x_i$ for all $i\in I\setminus\{k,h\}$. Note that in this case we must have $x'_h < x_h$, because $\emptyset\neq \{i\in I\mid x'_i<x_i\}\subseteq \{k,h\}$, $x_k=m(x)$, and $x'_h=m(x')$. Since $n\geq 3$, there is an index $j$ such that $j \in I\setminus \{k,h\}$.
Let us define positions $z,z'\in\ZZ_+^I$ by
\[
z_i = \begin{cases}
x_h+(x_0-x_0') & \mbox{if } i=h,\\
x_j+x_0' & \mbox{if } i=j,\\
x_i & \mbox{otherwise},
\end{cases}
~~~~
z_i' = \begin{cases}
x'_j+x_0' & \mbox{if } i=j,\\
x'_i & \mbox{otherwise}.
\end{cases}
\]
Then we have $z \in Z(x_0;x)$ and $z\to z'$ is a move in \NIM$_{\cH}$, because
\eqref{e-zpzxpx} holds in this case, too.

\smallskip
Finally note that in all four cases we have
$z' \in Z(x'_0;x')$ by the definitions of $z'$ and $Z(x'_0;x')$.
\qed

\bigskip

\noindent\emph{Proof of Theorem \ref{t1}}.
For a position $(x_0;x)$ of the extension,
let $B(x_0;x)=\{(x'_0;x') \mid (x_0;x)\to (x'_0;x')\}$.
Then we have
\[
\{(m(x'),u(x')+x_0') \mid (x'_0;x') \in B(x_0;x)\}=A(x_0;x),
\]
where $A(x_0;x)$ is as defined in \eqref{e-A}.
This is because $\{(m(x'),u(x')+x_0') \mid (x'_0;x') \in B(x_0;x)\}\supseteq A(x_0;x)$ follows by Lemmas
\ref{lemma-1az-1} and \ref{lemma-1az0}, and the opposite inclusion follows from Lemmas \ref{lemma-1az-1} and \ref{lemma-1az1}.

This completes the proof. \qed

\begin{example}
\label{example1}
Let us illustrate Theorem \ref{t1} and its Corollary \ref{c1} by some numerical examples for the extension $\Gamma_+$ of $\Gamma$,
the proper \NIM~ with $n=3$.

Consider first the position $(x_0;x)=(0;3,3,4)$. In this case we have $m(x)=3$, $y(0;x)=1$ and thus $z(0;x)=\binom{y(0;x)+1}{2}+1=2$.
Since $m\geq z$ we get for the SG function by \eqref{e-uv} that
\[\G_{\Gamma_+}(0;3,3,4)=(z-1)+((m-z)\mod (y+1))=1+(1\mod 2)=2.
\]
Note that for any move $x\to x'$ in $\Gamma_+$ (and since $x_0=0$, also in $\Gamma$) we have $u(x')\geq m(x)=3$.
Thus, to argue that \eqref{e-uv} provides the SG value of $2$ for this position,
we only need to consider moves to positions $x'$ for which $m(x')\geq z(0;x')$.
We list these positions (up to a permutation of the coordinates) in the table below:
\[
\begin{array}{c|c|c|c|l}
(0;x') & \,m(x')\, & \,y(0;x')\, & \,z(0;x')\, & \,\G=(z-1) + ((m-z) \mod (y+1)) \\
\hline
(0;3,2,2) & 2 & 1 & 2 & \,1=(2-1)+((2-2)\mod (1+1)) \\
(0;3,3,3) & 3 & 0 & 1 & \,0=(1-1)+((3-1)\mod (0+1))
\end{array}
\]
Let us next consider the position $(x_0;x_1,x_2,x_3)=(1;3,3,4)\in\ZZ_+^4$ of the extended game $\Gamma_+$.
In this case we have $m(x)=3$, $y(x_0;x)=2$, and thus $z(x_0;x)=\binom{2+1}{2}+1=4$.
Since $m(x)<z(x_0;x)$ we get by \eqref{e-uv} that
\[
\G_{\Gamma_+}(1;3,3,4)=u(x)+x_0=11.
\]
Note that for any move $(x_0;x)\to (x_0';x')$ we have
$
\G_{\Gamma_+}(x_0';x')\leq x_0'+u(x')<x_0+u(x)=11$.
We list below some moves $(x_0;x)\to (x_0';x')$ such that $m(x')<z(x_0';x')$ and the corresponding values by \eqref{e-uv}:
\[
\begin{array}{c|c|c|c|c}
(x_0';x') & \,m(x')\, & \,y(x_0';x')\, & \,z(x_0';x')\, & \,\G=u(x')+x_0' \\
\hline
(1;3,2,4) & 2 & 4 & 11 & 10 \\
(1;3,1,4) & 1 & 6 & 22 & 9\\
(1;3,0,4) & 0 & 8 & 37 & 8\\
(1;3,0,3) & 0 & 7 & 29 & 7\\
(1;3,0,2) & 0 & 6 & 22 & 6\\
(1;3,0,1) & 0 & 5 & 16 & 5\\
(1;3,0,0) & 0 & 4 & 11 & 4\\
(0;3,0,0) & 0 & 3 & 7 & 3
\end{array}
\]
Note finally that $(0;3,3,4)$ is reachable form $(1;3,3,4)$ and any position reachable form $(0;3,3,4)$ is also reachable form $(1;3,3,4)$, and thus the above computations show that
\[
11=\mex \{\G_{\Gamma_+}(x_0';x_1',x_2',x_3') \mid (1;3,3,4)\to (x_0';x_1',x_2',x_3')\}.
\]
\end{example}

\noindent\emph{Proof of Theorem \ref{tt2}}.
The same proof as above will work if we use $x_0=\T_\cF(x^J)$ and $x=x^I$. In particular, in the proof of Lemma \ref{lemma-1az0}
there exists a move $(x^J,x^I)\to ((x')^J,(x')^I)$ such that $\T_\cF((x')^J)=x'_0$ because $\cF$ is a Tetris hypergraph.
\qed


\section{Extension of $2$-pile \NIM}\label{s3}
Note that for $n=1$ the extension is equivalent with a $1$-pile \NIM.
In this section we consider the extension of proper \NIM~ for $n\leq 2$.
Surprisingly, this case seems to be much more difficult than the case of $n\geq 3$.
Here we present some partial results and conjectures.

\subsection{Upper and lower bounds for the SG function of extended proper \NIM}

For the analysis of the extension of $2$-pile \NIM,
we will need a few more properties of the SG function of extended proper \NIM.
Let us consider the hypergraph $\cH=2^I\setminus\{\emptyset,I\}$.

\begin{lemma}
\label{l-strong-mon}
SG function $\G_{\cH_+}(x_0;x)$ is
strictly monotone with respect to $x_0$.
\end{lemma}

\proof
Consider two positions  $(x_0;x)$  and
$(x'_0;x)$ such that $x_0 > x'_0$. Since $(x'_0;x)$ and any position reachable from
$(x'_0;x)$  are also reachable from $(x_0;x)$, we must have
$\G_{\cH_+}(x_0;x) > \G_{\cH_+}(x'_0;x)$ by the definition of the SG function.
\qed

To a position $(x_0;x)\in \ZZ_+^{I_+}$ let us associate
\[
\ell(x_0;x) = x_0 + \G_{\cH}(x).
\]

\begin{lemma}
\label{ul-bounds}
We have
\begin{equation}\label{e-lu}
\ell(x_0;x) \leq \G_{\cH_+}(x_0; x) \leq x_0+u(x).
\end{equation}
\end{lemma}

\proof
The upper bound is obvious, while
the lower one follows from Lemma \ref{l-strong-mon}.
\qed

In the rest of this section we consider $n=2$.
In this case the extended game is \NIM$_{\cH_+}$ where $\cH_+=\{\{0\},\{1\},\{2\},\{0,1\},\{0,2\}\}$. For simplicity, we change our notation. We denote by $x=(x_0,x_1,x_2)\in\ZZ_+^3$ a position of the extended game and by $\G(x)$ the value of its SG function.
Since proper \NIM~ with $n=2$ is the same as a $2$-pile \NIM,
we also use $u(x)=x_0+x_1+x_2$ and $\ell(x)=x_0+(x_1\oplus x_2)$, where $\oplus$ is the so called \NIM~ sum, see e.g., \cite{BCG01-04}.

In this case Lemma~\ref{ul-bounds} turns into the following inequalities

\begin{equation}\label{eq-ul-bounds}
\ell(x)  \leq \G(x)  \leq  u(x).
\end{equation}
If  $x_0 = 0$, the lower bound is attained.
Obviously, the upper bound is attained if  $x_1 = 0$  or  $x_2 = 0$.
We shall see that both bounds are attained in many other cases.


\subsection{Shifting  $x_1$  and  $x_2$  by the same power of  $2$} \label{ss21}

We present simple conditions under which the SG function
$\G(x)$  is invariant with respect
to a shift of $x_1$  and  $x_2$  by the same power of  $2$.

Let us note that $\G$ is not a monotone increasing function of $x_1$ and/or $x_2$, in contrast to $x_0$. The next lemma shows that a weaker property still holds, if we use power of $2$ increments.
\begin{lemma}
\label{lemma-2-a0}
For any
$k \in \Zp$, let us set  $\Delta^k=(0,2^k,2^k)$. Then,
\[
\begin{array}{rl@{\text{~~if~~~}}l}
\cG(x +\Delta^k) &= \cG(x) &  \cG(x) < 2^k, \text{ and}\\
\cG(x +\Delta^k) &\geq 2^k & \G(x)\geq 2^k.
\end{array}
\]
\end{lemma}

\proof
We show this by induction on $x$.
We first note that $\cG(0,0,0)=\cG(0,a,a)$ holds
for any positive integer $a$, which proves the base of the induction.

For a position $x$, we assume that the statement is true
for all $x'$ with $x' \leq x$, $x'\neq x$ and show that it holds
for $x$ by separately considering the cases  $\cG(x)<2^k$ and $\cG(x)\geq 2^k$.

Case 1: $\cG(x) < 2^k$.
We note that if $x \to x'$ is a  move then
so is $x+\Delta^k \to x'+\Delta^k$.
For any $v$ with $0 \leq v < \cG(x) < 2^k$, there exists
a move $x \to x'$ such that $\cG(x')=v$ by the definition of an SG function.
It follows from the induction hypothesis that $\cG(x'+\Delta^k)=\cG(x')=v$.
Since $x + \Delta^k \to x'+ \Delta^k$ is a  move, we have $\cG(x+\Delta^k)\neq v$.
Since this applies for values $0\leq v<\cG(x)$, we can conclude that
$\cG(x+\Delta^k)\geq \cG(x)$.

We will show next that for any position
$x'$  obtained from  $x + \Delta^k$  by a move
$x + \Delta^k \to x'$, the SG function values
of $x$ and $x'$ differ,  $\cG(x') \not= \cG(x)$.

Let $x'=(x'_0, x_1',x_2+2^k)$  without loss of generality.
If $x_1'< 2^k$, then we have
\[
\cG(x')\geq \ell(x')\geq x_1' \oplus (x_2+2^k) \geq 2^k > \cG(x).
\]
If  $x_1' \geq 2^k$ and $\cG(x'-\Delta^k)=\cG(x_0',x_1'-2^k,x_2) \geq  2^k$,
then  by induction hypothesis,
$\cG(x')=\cG((x'-\Delta^k)+\Delta^k)\geq 2^k$,
which implies that $\cG(x') \not=\cG(x)$.
Finally, if  $x_1' \geq 2^k$ and $\cG(x'-\Delta^k) < 2^k$,
then
$\cG(x'-\Delta^k)\not=\cG(x)$, since $x \to x'-\Delta^k$ is a  move.
By our induction hypothesis, we have $\cG(x')=\cG(x'-\Delta^k)\not=\cG(x)$.
This completes the proof of the case $\cG(x) < 2^k$.

Case 2: $\cG(x)\geq 2^k$.
By definition,  for any integer $v$ with $0 \leq v < 2^k \leq \cG(x)$, there exists
a  move $x \to x'$ such that $\cG(x')=v$.
By  induction hypothesis, we have  $\cG(x')=\cG(x'+\Delta^k)$.
Since $x'+\Delta^k$ is reachable from $x+\Delta^k$, we obtain $\cG(x+\Delta^k)\geq 2^k$.
\qed

To illustrate the first claim of the lemma, let
us note that:
\[
\begin{array}{rll}
\cG(1,0,0) &= 1 &< 2\\
\cG(2,0,0) &= 2 &< 4\\
\cG(3,0,0) &= 3 &< 4\\
\cG(4,0,0) &= 4 &< 8\\
\cG(0,1,2) &= 3 &< 4\\
\cG(1,0,1) &= 2 &< 4\\
\cG(1,1,1) &= 3 &< 4\\
\cG(1,1,2) &= 4 &< 8\\
\cG(2,1,3) &= 6 &< 8\\
\cG(3,1,1) &= 5 &< 8.
\end{array}
\]

Furthermore, computations show that
\[
\begin{array}{rll}
\cG(1,2,6) &= 5 &< 8 \\
\cG(1,5,6) &= 11 &< 16\\
\cG(2,5,5) &= 11 &< 16
\end{array}
\]
Hence,
for all
$i \in \Zp$  we have

\[
\begin{array}{rllll}
\cG(1,0,0) &= \cG(1,2,2) &= \cG(1,4,4) &= \cdots &= \cG(1,2i,2i) = 1\\
\cG(2,0,0) &= \cG(2,4,4) &= \cG(2,8,8) &= \cdots &= \cG(2,4i,4i) = 2\\
\cG(3,0,0) &= \cG(3,4,4) &= \cG(3,8,8) &= \cdots &= \cG(3,1+4i,1+4i) = 3\\
\cG(4,0,0) &= \cG(4,8,8) &= \cG(4,16,16) &= \cdots &= \cG(4,8i,8i) = 4\\
\cG(1,0,1) &= \cG(1,4,5) &= \cG(1,8,9) &= \cdots &= \cG(1,4i,1+4i) = 2\\
\cG(0,1,2) &= \cG(0,5,6) &= \cG(0,9,10) &= \cdots &= \cG(0,1+4i,2+4i) = 3\\
\cG(1,1,1) &= \cG(1,5,5) &= \cG(1,9,9) &= \cdots &= \cG(1,1+4i,1+4i) = 3\\
\cG(1,1,2) &= \cG(1,9,10) &= \cG(1,17,18) &= \cdots &= \cG(1,1+8i,2+8i) = 4\\
\cG(2,1,3) &= \cG(2,9,11) &= \cG(2,17,19) &= \cdots &= \cG(2,1+8i,3+8i) = 6\\
\cG(3,1,1) &= \cG(3,9,9) &= \cG(3,17,17) &= \cdots &= \cG(3,1+8i,1+8i) = 5\\
\cG(1,2,6) &= \cG(1,10,14) &= \cG(1,18,22) &= \cdots &= \cG(1,2+8i,6+8i) = 5\\
\cG(1,5,6) &= \cG(1,21,22) &= \cG(2,47,48) &= \cdots &= \cG(1,5+16i,6+16i) = 11\\
\cG(2,5,5) &= \cG(2,21,21) &= \cG(1,37,37) &= \cdots &= \cG(2,5 + 16i, 5 + 16i) = 11.
\end{array}
\]
\noindent

To illustrate the second claim of the lemma, let us
consider $k = 1$, $x = (1,2,3)$,  $x + (0,2,2) = (1,4,5)$,  and note that
$\cG(1,2,3) = 6$ while  $\cG(1,4,5)=2$.

\medskip

Lemma \ref{lemma-2-a0} results immediately the following claim.

\begin{corollary} \label{corollary-2-a1}
For any
$k \in \ZZ_+$,
if $\cG(x) < 2^k$ and $x \geq \Delta^k$ then $\cG(x-\Delta^k)=\cG(x)$.
\qed
\end{corollary}

We also need the following elementary arithmetic statement.

\begin{lemma} \label{lemma-2-a1}
Let $a$ be an integer with $2^{k-1} \leq a < 2^k$ for some positive integer $k$.
Then we have $a\oplus b>b$ if  $0\leq b < 2^{k-1}$,
and  $a\oplus b< b$ if  $ 2^{k-1} \leq b < 2^k$.
\end{lemma}

\proof
The claim results immediately from
the definition of the \NIM-sum, since
the binary representation of  $a$  includes $2^{k-1}$.
\qed

For a nonnegative integer $v$, let $k(v)$
denote the unique nonnegative integer such that  $2^{k(v)-1} \leq v < 2^{k(v)}$.

The following  consequence of Lemma \ref{lemma-2-a0} provides a characterization of the SG values.

\begin{theorem} \label{corollary-2-a2}
Let $x = (x_0, x_1, x_2)$  be a position with
$\cG(x) = v$  and  $x_1 \leq x_2$.
Then there exists a position  $\hat{x}$  such that
$\hat{x} \leq (v, 2^{k(v)-1}-1, 2^{k(v)}-1)$, $\cG(\hat{x})=\cG(x)$, and
$x=\hat{x}+\lambda \Delta^{k(v)}$ for some $\lambda\in \Zp$.
\end{theorem}

\proof
By Corollary \ref{corollary-2-a1}, there exists
a position $\hat{x}= (\hat{x}_0, \hat{x}_1, \hat{x}_2) = x-\lambda\Delta^{k(v)}$ for a nonnegative integer $\lambda$
such that $\cG(\hat{x})=\cG(x)=v$ and $\hat{x}_1< 2^{k(v)}$.
We show that $\hat{x}_0 \leq v$, $\hat{x}_1< 2^{k(v)-1}$, and
$\hat{x}_2< 2^{k(v)}$, which will imply the statement.
Since

\begin{equation}
\label{eq-2-ccc1}
v=\cG(\hat{x})\geq \ell(\hat{x})=\hat{x}_0+(\hat{x}_1 \oplus \hat{x}_2),
\end{equation}

\noindent
$ \hat{x}_0 \leq v$ holds.
Moreover, if $\hat{x}_2 \geq 2^{k(v)}$, then
$\hat{x}_1 \oplus \hat{x}_2 \geq 2^{k(v)}$ by $\hat{x}_1< 2^{k(v)}$,
which again contradicts \raf{eq-2-ccc1}, since $v < 2^{k(v)}$ by definition of $k(v)$.
Thus we have $\hat{x}_2 < 2^{k(v)}$.
Suppose that $ 2^{k(v)-1} \leq \hat{x}_1\leq \hat{x}_2$.
Then $\hat{x} \to (0,\hat{x}_1,\hat{x}_1\oplus v)$ is a  move, since
$\hat{x}_1\oplus v < \hat{x}_1 \leq \hat{x}_2$ by Lemma \ref{lemma-2-a1}.
This together with $\cG(0,\hat{x}_1,\hat{x}_1\oplus v)=v$ contradicts that $\cG(\hat{x})=v$.
\qed

\medskip

For any nonnegative integer $v$,  let us define
\[
Core(v) ~=~ \{x=(x_0, x_1,x_2)\mid \cG(x)=v, x_0\leq v, x_1< 2^{k(v)-1}, x_2< 2^{k(v)}\}.
\]
Then, Theorem \ref{corollary-2-a2} shows that
every position $x$ with $\cG(x)=v$ and $x_1 \leq x_2$ has a position $\hat{x}\in Core(v)$ such that
$x=\hat{x}+\lambda\Delta^{k(v)}$ for some nonnegative integer $\lambda$.

\bigskip

Note that $Core(v)$ has at most $2v^3$ positions, and by the definition of the SG function, we can compute their SG value in $O(v^5)$ time. This implies the following corollary.

\begin{corollary}
For any $v\in \Zp$ and for any position $x$ we can compute the value $\G(x)$, if $\G(x) \leq v$, or prove that $\G(x)>v$ in $O(v^5)$ time. \qed
\end{corollary}

This implies that we can compute the SG value of a position $x\in\ZZ^3$ in polynomial time in $\G(x)$, regardless the magnitude of the coordinates of $x$.


\subsection{Conjectures and partial results} \label{ss22}

In this subsection we assume that $x_1\leq x_2$.

\subsubsection{Case 1:  $x_1$ is a power of  $2$}
\label{PowerOf2}

Computational results suggest that
if  $x_1$  is a power of  $2$  then  $\G(x)$
equals either the lower or the upper bound in accordance with the following simple rule.

\begin{conjecture}
\label{t2}
Given  $x = (x_0, x_1, x_2)$  such that
$x_0 \in \ZZ_+$  and  $x_2 \geq x_1 = 2^k$  for some
nonnegative integer  $k$,
then  $\G(x) = \ell(x) = x_0 + (x_1 \oplus x_2)$ for any
\begin{equation}
\label{eq2^k}
x_2 = (2j+1) 2^k + m
\mbox{ such that }
j \in \ZZ_+
\mbox{ and }
0 \leq m < 2^k - x_0.
\end{equation}
Otherwise,
$\G(x) = u(x) = x_0 + x_1 + x_2$.
\end{conjecture}

Note that  \raf{eq2^k}  can be equivalently rewritten as
\[
2^k \leq (x_2 \mod {2^{k+1}}) < 2^{k+1} - x_0.
\]
It is also convenient to equivalently reformulate this conjecture
replacing   $\G(x)$  by  $\delta(x) = u(x) - \G(x)$.
For any nonnegative integer $a,b,c \in \Zp$  let us introduce the function
\[
f(a, b, c)  ~=~ \begin{cases}
1,  & \mbox{if } (c \mod a) \geq b,\\
0,  & \text{otherwise}.
\end{cases}
\]

In particular,  $f(a, b, c) = 0$  whenever  $b \geq a$.

It is not difficult to verify that Conjecture \ref{t2}
can be reformulated as follows:

Given  $x = (x_0, x_1, x_2)$  such that
$x_0 \in \ZZ_+$  and  $x_2 \geq x_1 = 2^k$  for some
$k \in \ZZ_+$, then
\begin{equation} \label{eq2^k-A}
\delta(x) =
2^{k+1} f(2^{k+1}, x_0 + 2^k, x_0 + x_2) = 2x_1 f(2x_1, x_0 + x_1, x_0 + x_2).
\end{equation}

We can illustrate this by several simple examples.
\[
\begin{array}{c|c|c|l}
x_0 & x_1 & x_2 & \hspace*{2cm}\delta(x)\\
\hline
0 & 1 & 1,2,3,... & 2,0,2,0,2,... \\
0 & 2 & 2,3,4,... & 4,4,0,0,4,4,0,0,4,4,...\\
0 & 4 & 4,5,6,... & 8,8,8,8,0,0,0,0,8,8,8,8,...\\
2 & 4 & 4,5,6,... & 8,8,0,0,0,0,0,0,8,8,0,0,0,0,0,0,...\\
\end{array}
\]

The upper bound is attained, whenever  $x_0 \geq x_1$; for instance
if  $x_0 = 2$  and  $x_1 \leq 2$  then  $\delta(x) \equiv 0$.

\bigskip

Note that Conjecture \ref{t2}, if true, would imply
the following useful addition to Lemma~\ref{lemma-2-a0}. Note first that for $x=(x_0,0,x_2)$ we have
$$\ell(x_0,0,x_2)= \G(x_0, 0 , x_2) = u(x_0, 0, x_2) = x_0 + x_2.$$
By the conjecture above we also have
$$\G(x_0, 2^k, x_2 + 2^k) = \ell(x_0, 2^k, x_2 + 2^k) = x_0+x_2,$$
whenever vector  $x' = (x_0, 2^k, x_2 + 2^k)$
satisfies condition \raf{eq2^k}. Otherwise,
$$\G(x_0, 2^k, x_2 + 2^k) = u(x_0, 2^k, x_2 + 2^k) = x_0 + x_2 + 2^{k+1}.$$

\noindent
The following examples illustrate the above statement:
By Lemma \ref{lemma-2-a0} we have the equalities
\[
\begin{array}{rlllll}
13 &= \G(0, 0, 13) &= \G(0, 16, 29) &= \G(0, 32, 45) &= \G(0, 64, 79) &=  \cdots\\
17 &= \G(1, 0, 16) &= \G(1, 32, 48) &= \G(1, 64, 80) &= \cdots \\
19 &= \G(2, 0, 17) &= \G(2, 32, 49) &= \G(2, 64, 81) &= \cdots
\end{array}
\]
In addition, our computations show the following equalities, in agreement with the above conjecture. For the lower bound equalities we have
\[
\begin{array}{rlllll}
13 &= \G(0, 0, 13) &=  \G(0, 2, 15) \\
17 &= \G(1, 0, 16) &=  \G(1, 2, 18) &= \G(1, 4, 20) &= \G(1, 8, 24) \\
19 &= \G(2, 0, 17) &=  \G(2, 4, 21) &= \G(2, 8, 25),
\end{array}
\]
while the upper bound is attained in the following cases
$\G(0, 1, 14) = 15$, $\G(0, 4, 17) = 21$, $\G(0, 8, 21) = 29$.
$\G(1, 1, 17) = 19$,  $\G(1, 16, 32) = 49$,
$\G(2, 1, 18) = 21$, $\G(2, 2, 19) = 23$, and  $\G(2, 16, 33) = 51$.

\subsubsection{Case 2:  $x_1$  is close to a power of  $2$}

Our computations indicate that
for a position $x = (x_0, x_1, x_2)$
the upper bound is attained
whenever the semi-closed interval
$(x_1, x_1 + x_0]$ contains a power of  $2$.
Let us recall that  $x_1 \leq x_2$  is assumed.

\begin{conjecture}
\label{conj1}
If  $x_1 < 2^k \leq x_0 + x_1$  for some  
$k \in \Zp$, then  $\G(x) = u(x)$.
\end{conjecture}

Instead, we are able to prove only the following special case.

\begin{proposition}
\label{prop-1}
If  $x_0 \geq 2^{k - 1}$  and  $x_1 < 2^k$  for
some  $k \in \Zp$, then
$\cG(x) =u(x)$.
\end{proposition}

\proof
We show the claim by induction on $u(x)$.
If $u(x) = 0$ (that is, $x=(0,0,0)$) then clearly
$\cG(x)=u(x)=0$.
Assuming that the claim holds for all $x$ with $u(x) \leq p-1$,
consider a position $x$ with $u(x)=p$.
We claim that for each integer $v$ with $x_0 \leq v < u(x)$,
there exists a  move $x \to x'$ such that $\cG(x')=v$.
This will imply  $\cG(x)=u(x)$, since
$x_0$ and $u(x)$ are lower and upper bounds for $\cG(x)$, respectively.

Let $x'=(x_0',x_1',x_2')$ be a position such that
$x'_0 = x_0$, $x_1' = x_1$, and $0 \leq x'_2 < x_2$.
Then $x \to x'$ is a  move such that
$x'$ satisfies $x'_0=x_0 \geq 2^{k-1}$ and $\min \{x_1',x_2'\}\leq x_1 < 2^k$.
Thus, by induction hypothesis, for any integer $v$
satisfying  $x_0 + x_1 \leq v < u(x)$
there exists a  move $x \to x'$  such that  $\cG(x')=v$.

Then, let us consider moves
$x \to x'$  such that $0 \leq x_0' < x_0$, $x_1'=x_1$, and $x_2'=0$.
By definition, we have $\cG(x')=u(x')$, which shows that
for any integer $v$ with $x_1 \leq v < x_0+x_1$,
there exists a  move $x \to x'$ such that $\cG(x')=v$.
Hence, our claim is proven if  $x_0 \geq x_1$, because $\ell(x) \geq x_0$.

If $x_0 < x_1$ then for each integer $v$ with $2^{k-1}  \leq v < x_1$
consider a position $x'$ such that $x'_0=0$,
$x_1'=x_1$, and $x_2'=x_1\oplus v$.
It follows from Lemma \ref{lemma-2-a1} that
$x_1\oplus v < x_1 \leq x_2$, which implies that $x'$ is reachable from $x$.
Since $\cG(x')=v$, and $\ell(x)\geq x_0\geq 2^{k-1}$, the proof is completed.
\qed

Note that this Proposition is a special case of Conjecture \ref{conj1} if we choose $k$ to be the smallest integer satisfying $x_0 \geq 2^{k - 1}$  and  $x_1 < 2^k$.

\begin{corollary}
If  $x_0 \geq x_1$ then $\cG(x)=u(x)$.
\qed
\end{corollary}


\subsubsection{Case 3:  $x_2$  is close to a multiple of a power of  $2$}
Let us summarize the previous results and conjectures:

\begin{itemize}
\item[(a)] If   $x_0 = 0$  then the lower bound is attained:  $\G(x) = x_1 \oplus x_2$.
\item[(b)] If  $x_1$  is a power of  $2$  then the condition of Conjecture \ref{t2} holds.
\item[(c)] If
$x_1 < 2^k \leq x_0 + x_1$  for some  $k \in \Zp$, then the condition of Conjecture \ref{conj1} holds.
\end{itemize}

Thus, we assume from now on that

\begin{equation}
\label{x_1<2^k}
2^{k-1} < x_1 < x_0 + x_1 < 2^k  \mbox{ for some } k \in \Zp.
\end{equation}

In this case, our computations show that
$\G(x) = u(x)$   whenever
$x_2$  differs from a multiple of $2^k$ by at most $x_0$.

\begin{conjecture}
\label{t4}
If  $x = (x_0, x_1, x_2)$  satisfies  \eqref{x_1<2^k}, $x_1 \leq x_2$, and
$$
j 2^k-x_0 \leq x_2 \leq j 2^k + x_0
$$
for some $j,k\in \ZZ_+$,
then $\G(x) = u(x)$.
\end{conjecture}

Note that assumption \eqref{x_1<2^k} is essential.
For example, by computations we have
$\G(1, 4, 20) = \ell(1, 4, 20) = 1 + (4 \oplus 20) = 17<u(1,4,20)$, while the other conditions of Conjecture \ref{t4} hold.


\subsubsection{Case 4: $x_2$  is large}
Although the SG function looks chaotic in general, it seems that the pattern becomes
much more regular when  $x_2$  is large enough.
Unfortunately, we cannot predict how large should it be or prove any observed pattern.

\begin{conjecture}
\label{cnj1o}
We have $\G(x) = u(x)$ whenever
\eqref{x_1<2^k} and the following two conditions hold simultaneously:
\begin{itemize}
\item [\rm{(i)}] either $x_0 > 1$,  or $x_0 = 1$ and $x_1$  is odd, and
\item [\rm{(ii)}] $x_2$ is sufficiently large.
\end{itemize}
\end{conjecture}

Let us start with several examples where $\G(x)=u(x)$ with  $x_0 = 1$  and odd  $x_1$, if $x_2> \tau(x_1)$ is large enough:
\[
\begin{array}{c|c|c}
x_1 & \tau(x_1) & \delta(1,x_1,\tau(x_1)) \\
\hline
5 &  14 & 1\\
9 &  94 & 1\\
11 & 30 & 1\\
13& 30 & 1\\
17& 446 & 1\\
19&158 & 1\\
21& 94 & 3\\
23&62 & 1\\
25& 126 & 1\\
27& 62 & 3\\
29& 30 & 10
\end{array}
\]
where $\delta(x)=u(x)-\G(x)$.
Note that we skip the values  $x_1 = 2^k - 1$, that is,  $x_1 = 1, 3, 7, 15, 31$ because those cases are covered by Conjecture \ref{conj1}. When $x_1$ is even and $x_0=1$ the computations show a chaotic behavior.

\medskip

It seems that for   $x_0 > 1$ the upper bound is attained sooner
(that is, for smaller  $x_2$)  and for both odd and even  $x_1$.
For  $x_1 < 2^5 = 32$  and  $x_0 = 2$  we obtain:
\[
\begin{array}{c|c|c}
x_1 & \tau(x_1) & \delta(1,x_1,\tau(x_1)) \\
\hline
5& 5 &1\\
 9& 45 &1\\
 10& 44 &1\\
 11& 13&1\\
 12&   13 &24\\
 13&   13 &2\\
 17&   125 &1\\
 18&   125&1\\
 19&   61  &1\\
 20&  93 &2\\
 21&   61 &1\\
 22&   61 &1\\
 23&  29 &2\\
 24&   92 &4\\
 25&    61 &1\\
 26&    61 &2\\
 27&   29 &1\\
 28&   29 &56\\
 29&   29 &4\\
\end{array}
\]
For  $x_0 = 3$  the upper bound is achieved even faster;
for  $x_1 < 2^5 = 32$  we have:
\[
\begin{array}{c|c|c}
x_1 & \tau(x_1) & \delta(1,x_1,\tau(x_1)) \\
\hline
9  &  28 &1\\
 10 &  20 &1\\
 11 & 12 &1\\
 12&   12 &24\\
17  &  92 &1\\
 18 &  92 &1\\
 19&   60 &1\\
 20   & 60 &1\\
 21  &  60 &1\\
 22 &  28 &2\\
 23  & 28 &2\\
24  &  56 &47\\
 25 &   28 &3\\
26 &   28 &3\\
 27 &    28 &3\\
 28 &   28 &56
\end{array}
\]
In all the above examples we skip the values of  $x_1$
such that  $x_1 < 2^k \leq x_1 + x_0$  for some  $k \in \Zp$, since
in this case the condition of Conjecture \ref{conj1} holds.
We also skip values of $x_1 = 2^k$ because the condition of
Conjecture \ref{t2} holds in this case.

Finally, let us consider the case when
$x_1$  is even and  $x_0 = 1$.

\begin{conjecture} \label{cnj1e}
Given  $x_0 = 1$  and an even  $x_1$  such that
$2^{k-1} < x_1 < 2^k$  for some  $k \in \Zp$, then
$\delta(x_0, x_1, x_2)$
takes only even values and
becomes periodic in $x_2$ with period  $2^k$,  when $x_2$  is large enough.
\end{conjecture}

The examples for  $k \leq 5$  are presented below.
Note that for  $x_1 = 2,4,8,16$  the condition in Conjecture \ref{t2} is satisfied and hence all remaining cases are presented in the table below.
For a threshold integer $\tau(x_1)$ we define
\[
\pi(x_1)=(\delta(1,x_1,i)\mid i=\tau(x_1)+1,\tau(x_1)+2,...),
\]
and write $(0,0,0,1,2,1,2)^*=((0)^3,(1,2)^2)^*$ for the infinite sequence $(0,0,0, 1, 2, 1, 2,0,0,0, 1, 2, 1, 2,...)$.

\medskip

\[
\begin{array}{c|c|c|c|c}
x_1 & \tau(x_1) & 2^k& \pi(x_1) &\delta(1,x_1,\tau(x_1)) \\
\hline
6&14&8&(0^3,4,(0,2)^2)^*&12\\
10&109&16&(4,(0^3,4)^2,0,12,(0,4)^2,0)^*&2\\
12&109&16&(2,0^5,6,8,2,(0,2,4,2,0,2,4)^*&2\\
14&30&16&(0^3,4,(0,2)^6)^*&28\\
18&446&32&((0^3,4)^4,(0,2)^8)^*&1\\
20&400&32&((0,2,0,6)^3,0,2,(0^5,6,0,2)^2,0,22)^*&3\\
22&94&32&((0^3,4,(0,2)^2)^2,(0,2)^8)^*&4\\
24&456&32&((0,2)^2,0,16,((0,2)^3,0,10)^2,0^9,16)^*&10\\
26&126&32&((0^3,4)^2,(0,2)^{12})^*&1\\
28&104&32&((0,2,0,6)^5,0,2,0^5,6,0,2,0,14)^*&3\\
30&30&32&(0^3,4,(0,2)^{14})^*&60
\end{array}
\]

Interestingly, 
$\delta(x) = 0$   whenever  $x_2$  is odd,
except for only one case:  $x_1 = 12$, when
$\delta(1, 12, x_2)$  takes non-zero values $8,4,4$  for
$x_2 = 5, 9, 13 \mod 16$, respectively.

\smallskip

Let us also recall that
for  $x_1 = 2^k$,  Conjecture \ref{t2} of  Section \ref{PowerOf2}
gives similar periodical sequences, which take only two values:
$\delta = 0$  and  $\delta = 2x_1$.

\bigskip

Before concluding the section, we remark that we cannot separately prove the conjectures stated above.
Indeed, to show that  $\G(x) = v$, we have to verify that for any nonnegative integer  $v' < v$, there exists a  move $x \to x'$ such that $\G(x') = v'$.
Thus to prove one of the conjectures  by induction, we may need all other conjectures.
It is natural to prove them simultaneously.
However, at this moment we cannot, since for example, we have no bounds for $\tau$ in Conjectures \ref{cnj1o} and  \ref{cnj1e}.

\section{Slow \NIM}
\label{s4}
In this section, we consider a variant of \NIM, so called slow \NIM.
A move in a hypergraph \NIM game {\sc Nim}$_\cH$ is called {\em slow}
if each pile is reduced by at most one token.
Let us restrict both players by their slow moves, then the obtained game
is called {\em slow hypergraph Nim}.
We study SG functions and loosing positions of slow Moore \NIM~ and slow exact \NIM,
where they respectively correspond to hypergraphs $\cH = \{H \subseteq I \mid 1 \leq |H| \leq k\}$ and
$\cH=\{H \subseteq I  \mid |H|=k\}$ for some $k \leq n$.
We provide closed formulas for  the SG functions of both games when  $n = k = 2$  and  $n = k+1 = 3$,
where we remak that the SG function for slow exact \NIM~ when $n = k = 2$ is trivial.
We also characterize loosing positions for slow Moore \NIM~ if either  $n \leq k+2$  or  $n = k+3 \leq 6$ holds.

Here we only present the results, where all the proofs  
can be found in the preprint by Gurvich and Ho \cite{BH15}.

\medskip

Given a position  $x = (x_1, \ldots, x_n) \in \ZZ_+^I$, we will
always assume that its coordinates are nondecreasing $x_1 \leq \cdots \leq x_n$.
The {\em parity vector} $p(x)$  is defined as the vector
$p(x) = (p(x_1), \ldots, p(x_n)) \in \{0,1\}^I$ such that
$p(x_i) = 0$  if  $x_i$  is even, and $p(x_i) = 1$  if  $x_i$  is odd.
It appears that the status of a position
$x$ in the slow Moore {\sc Nim} in the cases below is  defined by $p(x)$ .


\begin{proposition}\label{p7}

The SG function $\G$ for slow Moore \NIM~
for $n = k = 2$ and $n = 3, k = 2$  are uniquely defined by
$p(x)$  as follows:
\begin{enumerate}
\item [\rm{(i)}]   For $n = k = 2$,
    \begin{align*}
    \G(x) =
    \begin{cases}
    0, \text{ if } p(x)  = (0,0) \\
    1, \text{ if } p(x)  = (0,1) \\
    2, \text{ if } p(x)  = (1,1) \\
    3, \text{ if } p(x)  = (1,0).
    \end{cases}
    \end{align*}
\item [\rm{(ii)}]  For $n = 3$ and $k = 2$,
    \begin{align*}
    \G(x) =
    \begin{cases}
    0 \text{ if } p(x)  \in \{(0,0,0), (1,1,1)\} \\
    1 \text{ if } p(x)  \in \{(0,0,1), (1,1,0)\} \\
    2 \text{ if } p(x)  \in \{(0,1,1), (1,0,0)\}\\
    3 \text{ if } p(x)  \in \{(0,1,0), (1,0,1)\}.
    \end{cases}
    \end{align*}
\end{enumerate}
\end{proposition}

We next consider loosing positions of slow Moore \NIM.

\begin{proposition} \label{p8}
Consider a slow Moore \NIM~ when  $n \leq k+2$ or $n = k+3 \leq 6$.
Then for a position $x \in \ZZ_+^I$,
we have the following five cases.
\begin{enumerate} \itemsep0em
\item [$(1)$] for $n=k$, $x$ is loosing if and only if $p(x) = (0, 0, \ldots, 0)$.
\item [$(2)$] for $n=k+1$, $x$ is loosing  if and only if $p(x) \in \{(0, 0, \ldots, 0), (1, 1, \ldots, 1)\}$.
\item [$(3)$] for $n=k+2$, $x$ is loosing  if and only if $p(x) \in \{(0, 0, \ldots, 0), (0,1, \ldots, 1)$.
\item [$(4)$] for $n=5$ and $k=2$, $x$ is loosing  if and only if
    $$p(x) \in \{(0,0,0,0,0), (0,0,1,1,1), (1,1,0,0,1), (1,1,1,1,0)\};$$
\item [$(5)$] for $n=6$ and $k=3$, $x$ is loosing  if and only if
    $$p(x) \in \{(0,0,0,0,0,0), (0,0,1,1,1,1), (1,1,0,0,1,1), (1,1,1,1,0,0)\}.$$
\end{enumerate}
\end{proposition}

\noindent
We note  that Moore \NIM~ games satisfy that $1\leq k \leq n$, and
the case in which $n=4$ and $k=1$ is a standard $4$-pile \NIM.

\medskip

We finally consider slow exact \NIM.
Note that this game is trivial when $k=1$ or $k=n$.
We show that the game is tractable if $n=3$ and $k=2$.
Again the parity vector plays an important role,
although it does not define the SG function uniquely.

Define six sets of positions $x \in \ZZ_+^3$:
\begin{align*}
A   &= \{(2a, 2b-1, 2(b+i)) \mid     0 \leq a < b,    0 \leq i < a, (a+i)   \bmod{2} = 1\} \\
B   &= \{(2a,2b,2(b+i)+1) \mid       0 \leq a \leq b, 0 \leq i < a, (a+i)   \bmod{2} = 1\} \\
C_0 &= \{(2a-1, 2b-1, 2(b+i)-1) \mid 0 \leq a \leq b, 0 \leq i < a, (a+i)   \bmod{2} = 0\} \\
C_1 &= \{(2a-1, 2b-1, 2(b+i)-1) \mid 0 \leq a \leq b, 0 \leq i < a, (a+i)   \bmod{2} = 1\} \\
D_0 &= \{(2a-1, 2b, 2(b+i)) \mid   0 \leq a < b,    0 \leq i < a, (a+i) \bmod{2} = 1\} \\
D_1 &= \{(2a-1, 2b, 2(b+i)) \mid   0 \leq a < b,    0 \leq i < a, (a+i) \bmod{2} = 0\},
\end{align*}
and let $C   = C_1 \cup C_2$, $D  = D_1 \cup D_2$.

\begin{proposition} \label{p10}
For a slow exact \NIM~ with  $n=3$ and $k=2$,
the SG function $\G$ can be represented by
 \begin{align*}
    \G(x) =
    \begin{cases}
    0 \text{ if } x \in (\{(2a,2b,c) \mid 2a \leq 2b \leq c\} \setminus B) \cup A \cup C_0 \cup D_0 \\
    1 \text{ if } x\in (\{(2a,2b+1,c) \mid 2a \leq 2b+1 \leq c\} \setminus A) \cup B \cup C_1 \cup D_1 \\
    2 \text{ if } x \in \{(2a+1,2b+1,c) \mid 2a+1 \leq 2b+1 \leq c\} \setminus C \\
    3 \text{ if } x \in \{(2a+1,2b,c) \mid 2a+1 \leq 2b \leq c\} \setminus D.
    \end{cases}
    \end{align*}
\end{proposition}

\end{document}